\documentclass[a4paper]{article}
\usepackage[english]{babel}

\usepackage[authoryear,longnamesfirst]{natbib}
\usepackage[noblocks,auth-sc]{authblk}

\usepackage[margin=2.2cm]{geometry}
\usepackage[utf8]{inputenc}
\usepackage[T1]{fontenc}
\usepackage{amsmath,amsthm,amsopn,amsfonts,amssymb,dsfont,color}
\usepackage{mathrsfs}
\usepackage{tikz}
\usepackage{soul}

\usepackage{xcolor}
\usepackage{float}
\usepackage{hyperref}
\usepackage{dsfont}
\usepackage{stmaryrd}

\def\R {\mathbb R}

\def\H{\mathcal H}
\def\D{\mathcal D}

\def\Linf{L^\infty}
\def\m{M_0}
\def\moy{\langle f \rangle}
\def\vinf{v_\infty}
\def\uinf{u_\infty}

\def\HPhi{{\mathcal H}}

\def\R{\mathbb R}

\newcommand{\mo}{\vert \omega\vert}
\newcommand{\mO}{\vert \Omega\vert}

\newtheorem{theorem}{\textbf{Theorem}}[section]
\newtheorem{lemma}[theorem]{\textbf{Lemma}}

\theoremstyle{remark}

\numberwithin{equation}{section}

\title{\bf A functional inequalities approach for the  field-road diffusion model with (symmetric) nonlinear exchanges}
\date{}

\author[1,2]{Matthieu Alfaro}
\author[2]{Claire Chainais-Hillairet}
\author[3]{Flore Nabet}
\affil[1]{{\small {Univ. Rouen Normandie, CNRS, LMRS UMR 6085, F-76000 Rouen, France}}}
\affil[2]{{\small {Univ. Lille, CNRS, Inria, UMR 8524 - Laboratoire Paul Painlevé, F-59000 Lille, France}}}
\affil[3]{{\small {CMAP, CNRS, \'Ecole polytechnique, Institut Polytechnique de Paris, 91120, Palaiseau, France}}}

\begin{document}

\maketitle
\vspace{3pt}

\begin{abstract} 
In this note, we  consider the so-called  field-road  diffusion model in a bounded domain, consisting of  two parabolic PDEs posed on sets of different dimensions and  coupled through (symmetric) {\it nonlinear} exchange terms.  We propose a new and rather direct functional inequalities approach to prove the exponential decay of a relative entropy, and thus the convergence of the solution towards the stationary state selected by the total mass of the initial datum. 

\end{abstract}

\section{Introduction}\label{s:intro}

In this note, we are interested in the large time behavior of the solution  $(v,u)=(v(t,x,y),u(t,x))$  to the so-called {\it field-road diffusion model}
\begin{subequations}\label{syst}
\begin{align}
	&\partial_{t} v = d \Delta v,
&&\quad t>0, \; x \in \omega, \; y\in(0,L), \label{syst-eqv}\\
	&- d\, \partial_{y} v|_{y=0}= \alpha(\mu_0 u^\beta - \nu_0 (v|_{y=0})^\alpha),
&&\quad t>0, \; x \in \omega,\label{syst-CLv}\\
	&\partial_{t} u = D \Delta u + \beta ( \nu_0 (v|_{y=0})^\alpha - \mu_0 u^\beta),
&&\quad t>0, \; x \in \omega, \label{syst-equ}\\
 &\frac{\partial u}{\partial n'} = 0,
&&\quad t>0, \; x \in \partial \omega,\label{syst-CLu}\\
 &\frac{\partial v}{\partial n} = 0,
&&\quad t>0, \; x \in \partial \omega, \;y\in(0,L), \text{ and } \; x\in \omega, \; y=L,
\label{syst-CLv2}\end{align}
\end{subequations}
supplemented with an initial condition $(v_0,u_0) \in \Linf(\Omega)\times \Linf(\omega)$. Here, $\Omega\subset \R^N$ ($N\geq 2$)  is a bounded cylinder (the field) of the form
$$
\Omega=\omega \times (0,L), \quad \omega \text{ a bounded convex and open set    of  $\R^{N-1}$ (the road)},\; L>0.
$$
The unknowns $v$ and $u$ correspond to the densities of individuals, respectively in the field $\Omega$ and on the road $\omega$ ; $d$ and  $D$ are the (positive) diffusion coefficients in the field and on the road. Obviously, $\Delta v$ has to be understood as $\Delta_x v+\partial_{yy}v$, while $\Delta u$ has to be understood as $\Delta_x u$. For $u$ we impose the zero Neumann boundary conditions on the boundary $\partial \omega$ ($n'$ denotes the unit outward normal vector to $\partial \omega$). For $v$, we impose the zero Neumann boundary conditions on the lateral boundary $\partial \omega \times (0,L)$ and on the upper boundary $\omega \times \{L\}$ ($n$ denotes the unit outward normal vector to $\partial \Omega$). On the lower boundary $\omega \times \{0\}$, the exchanges between the field and the road correspond to the value of the  outward flux of $v$ given by \eqref{syst-CLv}, where $\mu_0>0$ and $\nu_0 >0$ are transfer coefficients. These exchange terms also imply the zeroth-order term in \eqref{syst-equ}, linking the field and the road equations ;  they are the core of the model and are here assumed to be {\it nonlinear}, namely $(\alpha,\beta)\in [1,+\infty)^2\setminus\{(1,1)\}$. 

\medskip

The field-road model was introduced  by 
\cite{Ber-Roq-Ros-13-2, Ber-Roq-Ros-13-1, Ber-Roq-Ros-shape, Ber-Roq-Ros-tw} as a model for the spreading of  diseases or invasive species in presence of networks with fast propagation (typically $D>d$). We refer to the introduction in \cite{Alf-Cha-25} for more details and references. 

Very recently, a series of works has focused on the  {\it purely diffusive} field-road system: the fundamental solution was obtained in  \cite{Alf-Duc-Tre-23},  the PDE model was retrieved as the hydrodynamic limit of a particle system in \cite{Alf-Mou-Tre-25}. The system \eqref{syst} with linear exchanges ($\alpha=\beta=1$) was studied in 
\cite{Alf-Cha-25}: in both the continuous and  the discrete settings, 
the dissipation of a quadratic entropy is proved, the main tool being an adapted  Poincaré-Wirtinger inequality, see Lemma \ref{lem:PW} below.

Note that \eqref{syst} stands in the class of  {\it volume-surface} systems, considered in \cite{Fel-Lat-Tan-18}, \cite{Egg-et-al-18}, and the references therein. In \cite{Fel-Lat-Tan-18}, nonlinear exchanges are considered and the long time behavior is studied thanks to the  (logarithmic) Boltzmann entropy. The proof is rather lengthy and intricate and our goal here is to provide a  more direct approach.

To do so, in this note, we consider the case of {\it symmetric} nonlinear exchanges, namely
\begin{equation}\label{alpha-egal-beta}
    \alpha=\beta>1.
\end{equation}
For simplicity we use the shortcuts $\mu=\alpha \mu_0$, $\nu=\alpha\nu _0$ so that the right hand side of the  equation \eqref{syst-CLv} is reduced to $\mu u^\alpha-\nu (v|_{y=0})^\alpha$. In this framework, we consider a (power) Tsallis  entropy  and rely on two distinct functional inequalities: a Poincaré-Wirtinger inequality adapted to the field-road coming from \cite{Alf-Cha-25} and a Beckner type inequality \cite{Bec-89} coming from \cite{Cha-Jun-Sch-16}. We believe the arguments become transparent and, furthermore, this approach can be transferred to the design of a numerical scheme preserving the main properties of the system, see \cite{Alf-Cha-25} for the linear case.

We hope the method to be adapted to even more complex situations. In particular we aim at addressing the issue of {\it nonsymmetric} exchanges, in the sense that $\alpha\neq \beta$, considered in \cite{Fel-Lat-Tan-18}. This still requires an improvement of our method.

\medskip

This note is organized as follows. In Section \ref{s:setting} we present the result, which is proved in Section \ref{s:proof} and completed by numerical explorations in Section \ref{s:num}.

\section{Setting of the result}\label{s:setting}

We start with some basic facts. We consider  $v_0\in L^\infty(\Omega)$, $u_0\in L^\infty(\omega)$, both nonnegative and not simultaneously trivial. As a result, the total mass is initially positive
\begin{equation*}
\m:=\int _\Omega v_0(x,y)\,dxdy+\int_\omega u_0(x)\,dx >0.
\end{equation*}
The definition, existence and uniqueness of the weak solution to the Cauchy problem follows the same lines as in \cite{Fel-Lat-Tan-18}, see also  \cite{Egg-et-al-18}, \cite{Alf-Cha-25}. We denote $(v=v(t,x,y), u=u(t,x))$ the weak  solution  starting from $(v_0=v_0(x,y), u_0=u_0(x))$. Since the initial data are nonnegative and bounded, it follows from the comparison principle that both $v$ and $u$ are nonnegative and uniformly bounded (by some constant depending on $\Vert v_0\Vert _{L^\infty}$, $\Vert u_0\Vert_{L^\infty}$, $\mu$, $\nu$, $\alpha$). The total mass of the system  $\int_\Omega v(t,x,y)\,dxdy+\int_\omega u(t,x)\,dx$ is constant, namely
\begin{equation}\label{masse-conservee}
\int_\Omega v(t,x,y)\,dxdy+\int_\omega u(t,x)\,dx=M_0, \quad \forall t>0.
\end{equation}
The unique constant steady-state $(\vinf,\uinf)$ with mass $\m$ is given by
\begin{equation}\label{eq-pour-steady}
\nu \vinf^\alpha=\mu \uinf^\alpha, \quad \mO \vinf+\mo \uinf=\m.
\end{equation}

We apply a relative entropy method as presented for instance in the book by  \cite{Jungel_2016}. Define 
$$\Phi(s):=\frac{s^{\alpha+1}-(\alpha+1)s}{\alpha}+1,
$$
which satisfies $
\Phi''>0$ , $\Phi'(1)=0$, $\Phi(1)=0$. We define a nonnegative entropy, relative to the steady-state $(\vinf,\uinf)$, by
\begin{equation}
\label{def-H}
\HPhi(t):=\int _\Omega \vinf \Phi\left(\frac{v(t,x,y)}{\vinf}\right)\,dxdy+\int_\omega \uinf \Phi \left(\frac{u(t,x)}{\uinf}\right)\,dx,
\end{equation}
which, using \eqref{masse-conservee} and \eqref{eq-pour-steady}, can be recast
\begin{equation}
\label{H-recast}
\HPhi(t)=\frac{1}{\alpha \vinf ^\alpha}\int _\Omega \left(v^{\alpha+1}(x,y)-\vinf^{\alpha+1}\right)\,dxdy+\frac{1}{\alpha \uinf ^\alpha}\int_\omega \left(u^{\alpha+1}(x)-\uinf^{\alpha+1}\right)\,dx.
\end{equation}

The main result then writes as follows. 

\begin{theorem}[Exponential decay of entropy]\label{th:main} 
Assume $\alpha=\beta>1$. Let $v_0\in L^\infty(\Omega)$ and $u_0\in L^\infty(\omega)$ be both nonnegative and not simultaneously trivial. Let $(v=v(t,x,y), u=u(t,x))$ be the  solution  to \eqref{syst} starting from $(v_0=v_0(x,y), u_0=u_0(x))$,  and $(\vinf,\uinf)$ the associated steady-state defined by \eqref{eq-pour-steady}.  Then the entropy defined by \eqref{def-H} decays exponentially, namely
\begin{equation}\label{H-decay}
0\leq \H (t)\leq \H (0) e^{-\lambda t}, \quad \forall t> 0,
\end{equation}
for some positive $\lambda=\lambda(N,\Omega,\mu,\nu,d,D,v_0,u_0,\alpha)$. 

As a by-product, 
\begin{equation}
\label{by-product}
\Vert v-\vinf \Vert_{L^{\alpha+1}(\Omega)}^{\alpha+1}+\Vert u-\uinf \Vert_{L^{\alpha+1}(\omega)}^{\alpha+1}\leq  M e^{-\lambda t}, \quad \forall t> 0,
\end{equation}
for some positive $M=M(N,\Omega,\mu,\nu,d,D,v_0,u_0,\alpha)$.
\end{theorem}

\section{A functional inequalities approach}\label{s:proof}

We prove Theorem \ref{th:main}. In the sequel, the notation $\mathcal A\lesssim \mathcal B$ means that $\mathcal A\leq C\mathcal B$ for some positive constant $C=C(N,\Omega,\omega,\mu,\nu,d,D,v_0,u_0,\alpha)$.

\medskip

By differentiating expression \eqref{H-recast} with respect to time, using the equations in \eqref{syst} and integration by parts we reach
\begin{multline*}
\frac{\alpha}{\alpha+1} \frac{d}{dt} \HPhi (t)=\frac{1}{\vinf^\alpha}\left(\int_\omega (\mu u^\alpha-\nu (v|_{y=0})^\alpha)v^{\alpha}-d\int_\Omega \alpha \nabla v\cdot v^{\alpha-1}\nabla v\right)\\
+\frac{1}{\uinf^\alpha}\left(-D\int_\omega \alpha \nabla u\cdot u^{\alpha-1}\nabla u +\int_\omega (\nu (v|_{y=0})^\alpha-\mu u^\alpha)u^{\alpha}\right).
\end{multline*}
Thanks to \eqref{eq-pour-steady} we can gather the two non gradient terms and obtain
\begin{equation}
\label{dHdt}
 \frac{d}{dt} \HPhi(t) \lesssim -\int _\Omega \vert \nabla (v^{\frac{\alpha+1}{2}})\vert ^2 -\int _\omega \vert \nabla (u^{\frac{\alpha+1}{2}})\vert ^2-\int_\omega \left(\nu (v|_{y=0})^\alpha-\mu u^\alpha\right)^2=:-\mathcal D(t).
\end{equation}

We now take advantage of the adapted Poincaré-Wirtinger inequality developed in \cite{Alf-Cha-25}. To do so, for $\ell>0$, we \lq\lq enlarge'' $\Omega=\omega \times (0,L)$ to $\Omega^+=\omega\times (-\ell,L)$. We denote $\Omega_\ell =\omega \times (-\ell,0)$ the so-called thickened road. We work with 
\begin{equation}
\label{d-sigma}
d\rho= \left(
\frac{\vinf}{\m}\mathbf 1_{\Omega}(x,y)+\frac{1}{\ell}\frac{\uinf}{\m}\mathbf 1_{\Omega_\ell}(x,y)
\right)\,dxdy,
\end{equation}
which is a probability measure as can  be checked thanks to \eqref{eq-pour-steady}, and with
\begin{equation}
\label{f(x,y)}
f(x,y)=\left(\frac{v(x,y)}{\vinf}\right)^\alpha \mathbf 1_{\Omega}(x,y)+ \left(\frac{u(x)}{\uinf}\right)^\alpha \mathbf 1_{\Omega_\ell}(x,y), \quad (x,y)\in \Omega^+,
\end{equation}
where we have omitted to write the $t$ variable.  By strictly reproducing the proof of \cite[Theorem 1]{Alf-Cha-25}, which was concerned with the case $\alpha=1$, we obtain the following.

\begin{lemma}[Adapted Poincaré-Wirtinger inequality]\label{lem:PW}
Defining $
 \moy:=\int_{\Omega ^+} f\, d\rho$, there holds
\begin{equation}\label{PW}
\Vert f-\moy\Vert_{L^2(\Omega^+,d\rho)} ^2 \lesssim \int_\Omega \vert \nabla (v^\alpha)\vert ^2+\int_\omega \vert \nabla (u^\alpha)\vert^2+ \int _\omega \left( \nu (v|_{y=0})^\alpha-\mu u^\alpha \right)^2.\end{equation}
\end{lemma}

Next, we borrow \cite[Lemma 7]{Cha-Jun-Sch-16}.

\begin{lemma}[Generalized Beckner inequality II]\label{lem:Beckner} For $0<q<2$, $pq\geq 1$, there holds
\begin{equation}
\label{BecII}
\Vert f\Vert _{L^q((\Omega^+,d\rho)}^{2-q}\left( \int _{\Omega^+}\vert f\vert ^q d\rho-\left(\int _{\Omega^+}\vert f\vert ^{\frac 1p} d\rho\right)^{pq}\right)\lesssim \Vert f-\moy\Vert_{L^2(\Omega^+,d\rho)} ^2.
\end{equation}
\end{lemma}

Now, observe that since $v$ and $u$ are uniformly bounded there holds 
$$
\vert \nabla (u^\alpha)\vert\lesssim \vert \nabla (u^{\frac{\alpha+1}{2}})\vert, \quad \vert \nabla (v^\alpha)\vert\lesssim \vert \nabla (v^{\frac{\alpha+1}{2}})\vert.
$$
It therefore follows from Lemma \ref{lem:PW} and Lemma \ref{lem:Beckner} (with $q=\frac{\alpha+1}{\alpha}$, $p=\alpha$) that
$$
\mathcal D(t)
\gtrsim  \Vert f\Vert _{L^{\frac{\alpha+1}{\alpha}}(\Omega^+,d\rho)}^{\frac{\alpha-1}{\alpha}}\left( \int _{\Omega^+}\vert f\vert ^{\frac{\alpha+1}{\alpha}} d\rho-\left(\int _{\Omega^+}\vert f\vert ^{\frac 1{\alpha}} d\rho\right)^{\alpha+1}\right).
$$
Now using \eqref{f(x,y)}, \eqref{d-sigma}, and \eqref{masse-conservee}, we see that
\begin{eqnarray*}
\int _{\Omega^+}\vert f\vert ^{\frac{\alpha+1}{\alpha}} d\rho-\left(\int _{\Omega^+}\vert f\vert ^{\frac 1{\alpha}} d\rho\right)^{\alpha+1}&=&\int _\Omega \left(\frac{v}{\vinf}\right)^{\alpha+1}\frac{\vinf}{M_0}dxdy+\int _\omega \left(\frac{u}{\uinf}\right)^{\alpha+1}
\frac{\uinf}{M_0}dx-1\\
&=& \frac{1}{M_0 \vinf ^\alpha}\int _\Omega \left(v^{\alpha+1}-\vinf^{\alpha+1}\right)+\frac{1}{M_0 \uinf ^\alpha}\int_\omega \left(u^{\alpha+1}-\uinf^{\alpha+1}\right),
\end{eqnarray*}
by using the second relation in \eqref{eq-pour-steady}. Next, by Jensen's inequality,
$$
\Vert f\Vert _{L^{\frac{\alpha+1}{\alpha}}(\Omega^+,d\rho)}^{\frac{1}{\alpha}}=\left(\int_{\Omega ^+} f^{\frac{\alpha+1}{\alpha}}d\rho \right)^{\frac 1{\alpha+1}}\geq \int_{\Omega ^+} f^{\frac 1{\alpha}} d\rho=1.
$$
As a result, we end up with
\begin{equation}\label{we-end-up}
\D(t)\gtrsim \frac{1}{M_0 \vinf ^\alpha}\int _\Omega \left(v^{\alpha+1}-\vinf^{\alpha+1}\right)+\frac{1}{M_0 \uinf ^\alpha}\int_\omega \left(u^{\alpha+1}-\uinf^{\alpha+1}\right)=\frac{\alpha}{M_0}\HPhi(t).
\end{equation}
In view of \eqref{we-end-up} and \eqref{dHdt}, we collect $ \frac{d}{dt} \HPhi(t) \lesssim -\HPhi(t)$, which proves \eqref{H-decay}. 

Last  the decay of the $L^{\alpha+1}$ norm, namely  \eqref{by-product}, follows from \eqref{def-H} and the fact that $\vert 1-s\vert ^{\alpha+1}\leq s^{\alpha+1}-(\alpha+1)s+\alpha=\alpha \Phi(s)$ for all $s\geq 0$.

Theorem \ref{th:main} is proved. \qed 

\section{Numerical experiments}\label{s:num}

Our aim in this section is to illustrate the exponential decay of the relative entropy ${\mathcal H}$ defined by \eqref{def-H}, as stated in Theorem \ref{th:main} for the nonlinear 
field-road model \eqref{syst} with symmetric exchanges ($\alpha=\beta$). We will also investigate the behaviour of the similar relative entropy in the case with nonsymmetric exchanges 
($\alpha\neq \beta$). In order to do some numerical investigations, we use a two-point flux approximation (TPFA) finite volume scheme, with a backward in time Euler method, as introduced in \cite {Alf-Cha-25} 
for $\alpha=\beta =1$. Due to the nonlinear exchanges, the scheme consists in a nonlinear system of equations at each time step, which is solved using Newton's method. 

For the numerical experiments, we consider that the one-dimensional road is $\omega=(-2L,2L)$ and the two-dimensional field  is $\Omega=\omega\times (0,L)$ with $L=20$. 
The value of the kinetical parameters $\mu_0$ and $\nu_0$ are $\mu_0= 1$, $\nu_0 = 5$. The diffusion parameters, in the field and in the road, are respectively 
$d=1$, $D=1$. We consider two test cases already proposed in \cite{Alf-Cha-25} and defined in Table \ref{table:cas-tests}. In both test cases, the individuals are scattered in the field, 
but the road is empty in Test case 1, while there are some individuals scattered in the road in Test case 2.

\begin{table}[htbp!]
\centering
\caption{Presentation of the test cases used for the numerical experiments.}
\begin{tabular}{|c|c|c|}
\hline
& Test case 1&Test case 2\\
\hline
$v_0(x,y)$ & {\scriptsize$\phantom{\dfrac{1}{2}}100\cdot {\mathbf 1}_{[-10,-7.5]\cup[-5,-2.5]\cup[2.5,5]\cup[7.5,10]}(x)\cdot {\mathbf 1}_{[7.5,10]}(y)$ } & {\scriptsize$150\cdot{\mathbf 1}_{[-10,-7.5]\cup[-5,-2.5]\cup[2.5,5]\cup[7.5,10]}(x)\cdot {\mathbf 1}_{[8.75,10]}(y) $ }\\
\hline
$u_0(x)$& 0 & {\scriptsize$\phantom{\dfrac{1}{2}}62.5\cdot {\mathbf 1}_{[-10,-7.5]\cup[-5,-2.5]\cup[2.5,5]\cup[7.5,10]}(x)$ } \\
\hline
\end{tabular}
\label{table:cas-tests}
\end{table} 

\begin{figure}[htb!]
\begin{tabular}{cc}
\includegraphics[width=0.45\textwidth]{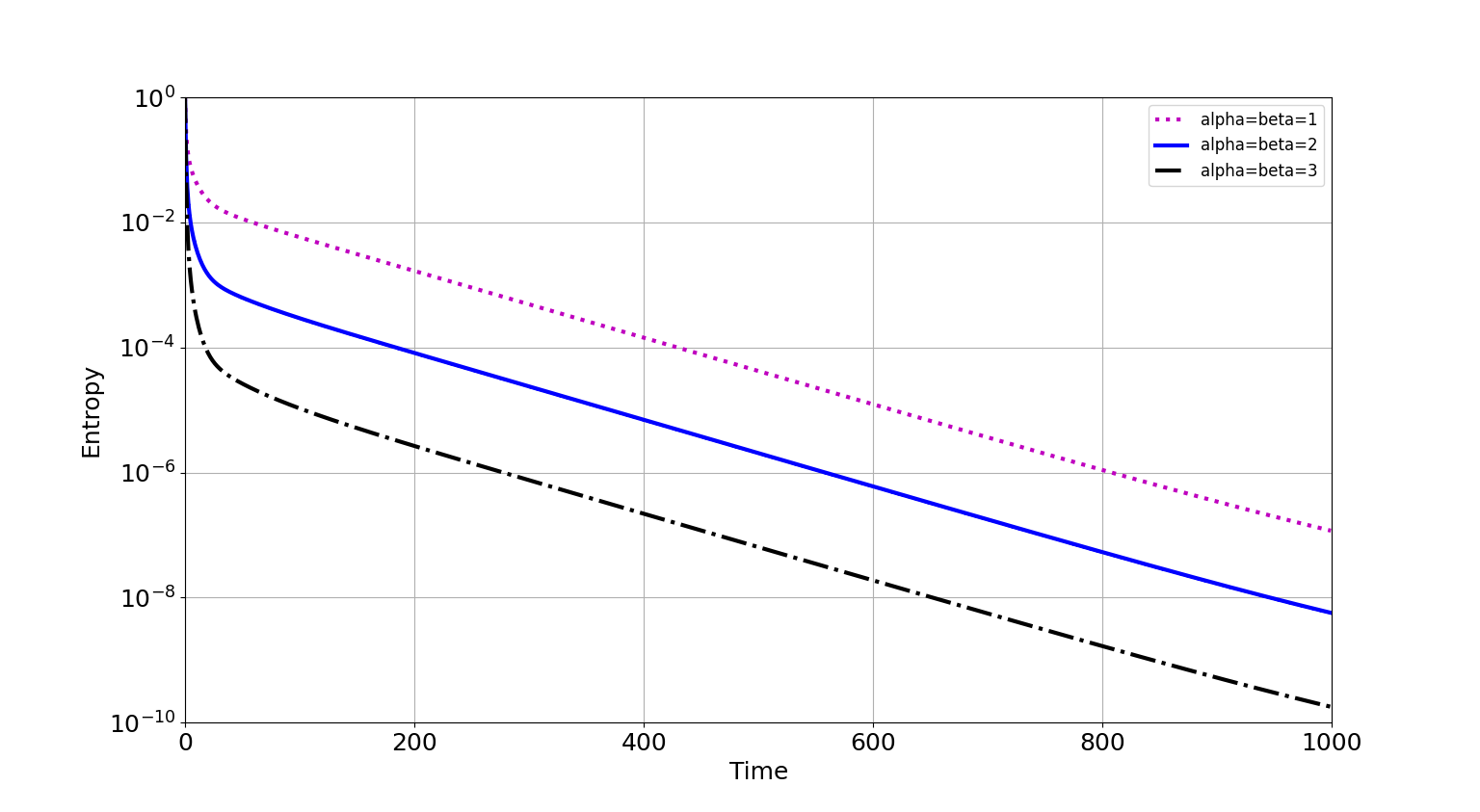}&
\includegraphics[width=0.45\textwidth]{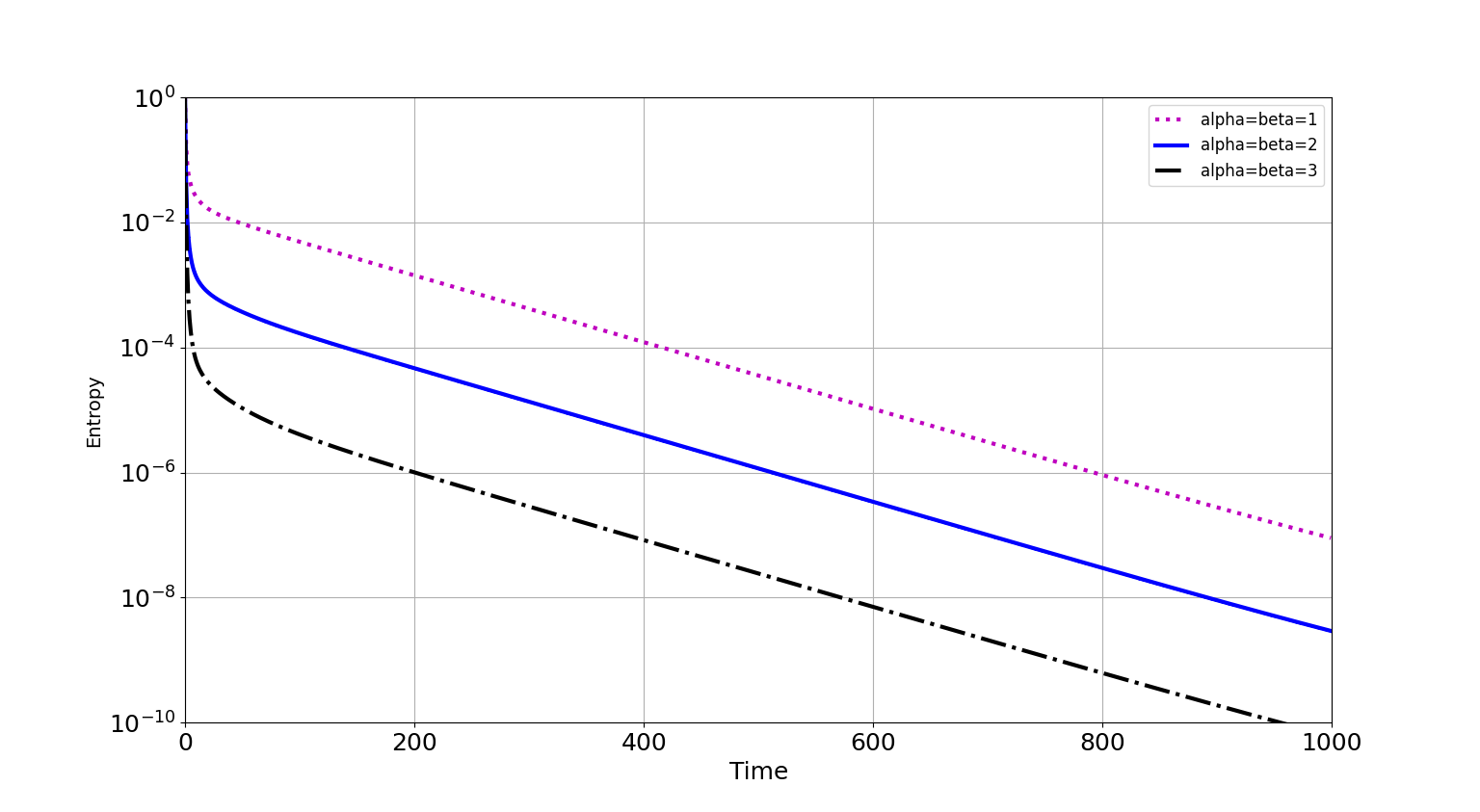}
\end{tabular}
\caption{Exponential decay of the relative entropy  in the symmetric case, $\alpha=\beta$, for Test case 1 (left) and Test case 2 (right). }\label{fig:sym}
\end{figure}

\begin{figure}[htb!]
\begin{tabular}{cc}
\includegraphics[width=0.45\textwidth]{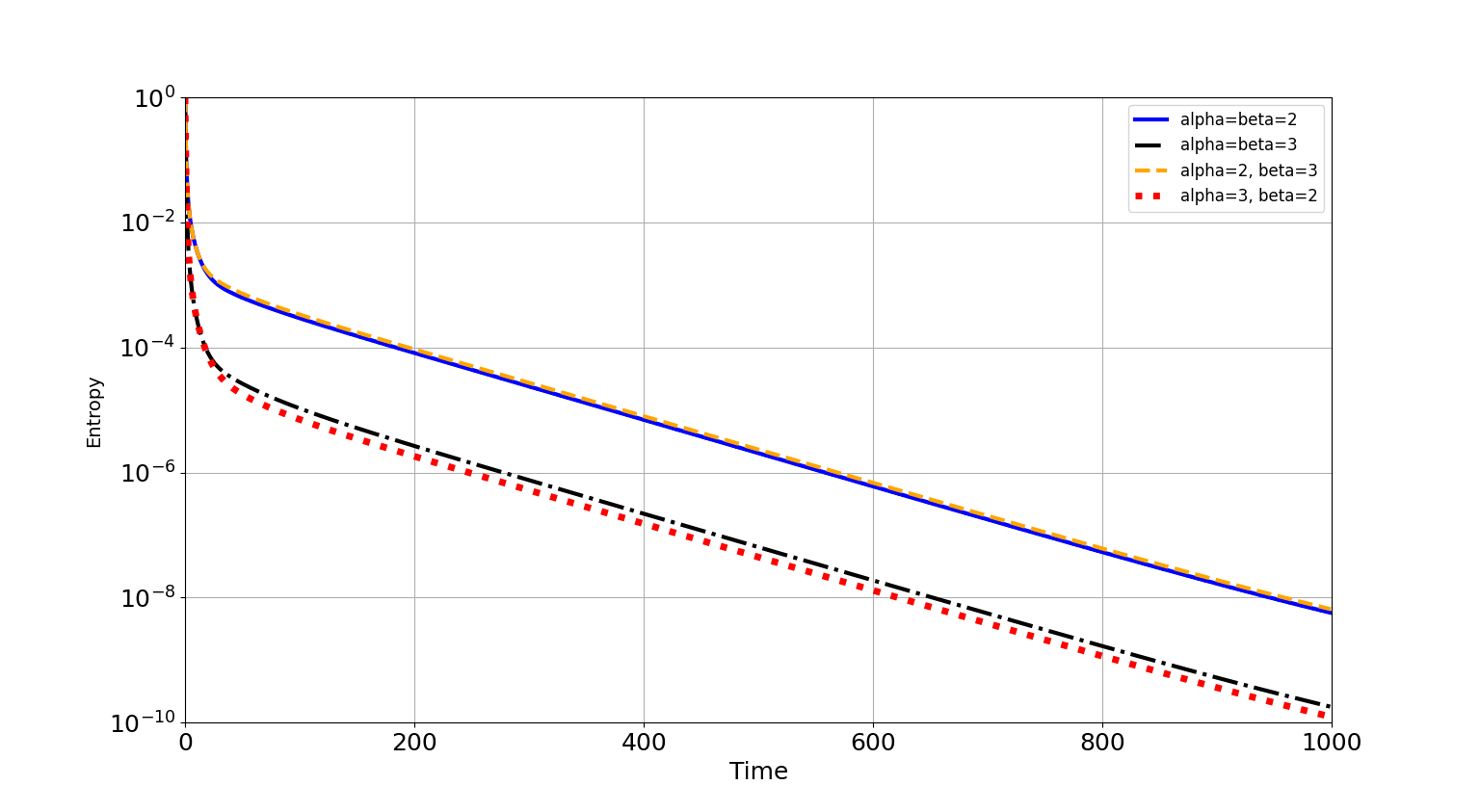}&
\includegraphics[width=0.45\textwidth]{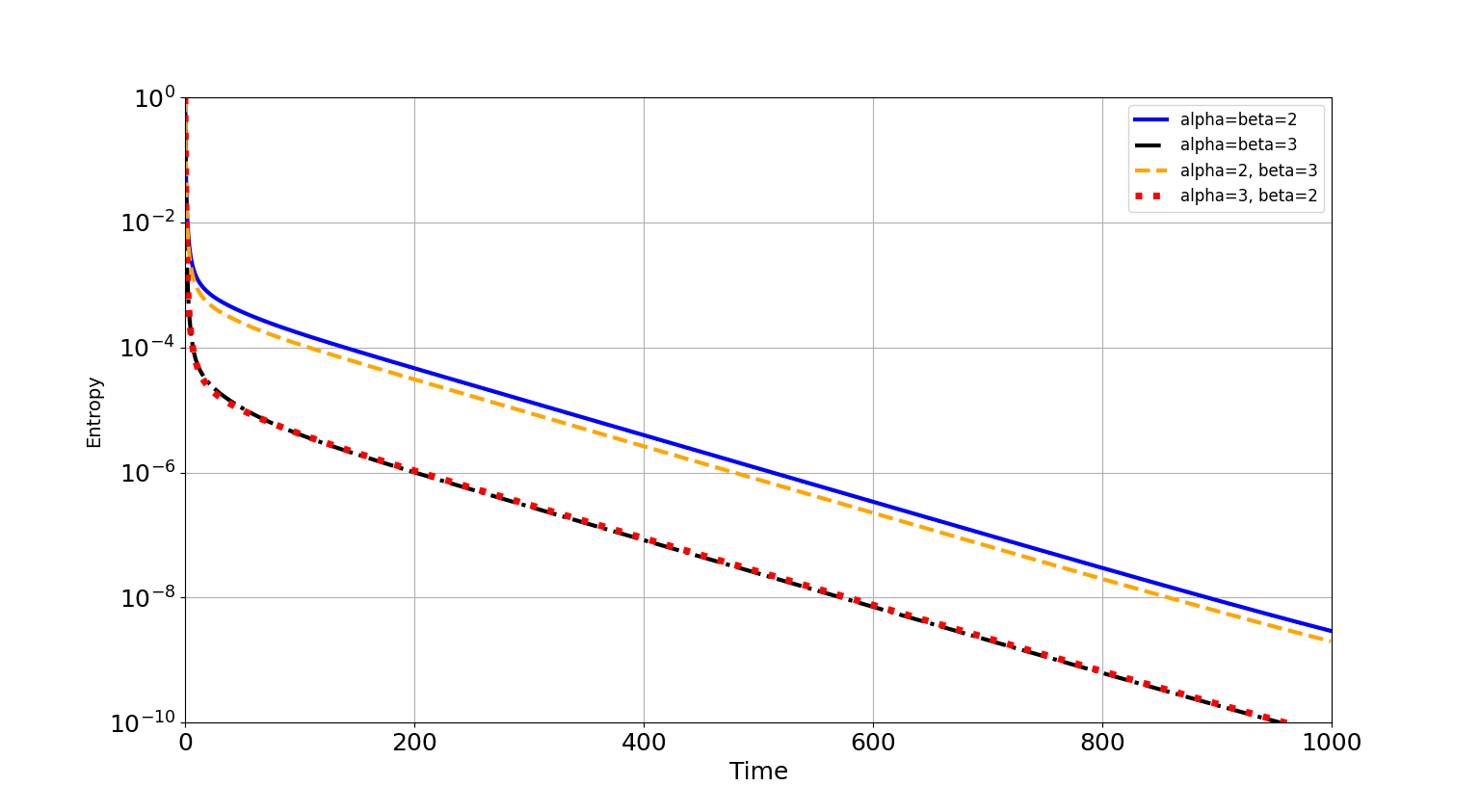}
\end{tabular}
\caption{Comparison of the decay of the relative entropy  in the symmetric  and nonsymmetric cases, for Test case 1 (left) and Test case 2 (right). }\label{fig:nonsym}
\end{figure}

Figure \ref{fig:sym} shows the long-time behaviour of the relative entropy ${\mathcal H}$  for both test cases in the symmetric case. Figure~\ref{fig:nonsym} shows the same 
evolution but in the nonsymmetric case. In the symmetric as in the nonsymmetric case, we observe that the decay of the relative entropy in time is exponential. 
Moreover, we observe that the decay rate seems to be independent of the values of $\alpha$ and $\beta$. It is also independent of the initial condition.

\bibliographystyle{elsarticle-harv}

\bibliography{biblio}

\end{document}